%%propp.tex: 
%%a Plain TeX file by Shalosh B. Ekhad and Doron Zeilberger (x pages)

%begin macros

\baselineskip=14pt
\parskip=10pt

\magnification=\magstephalf

\def\1{{\overline{1}}}
\def\2{{\overline{2}}}
\parindent=0pt
\overfullrule=0in

\def\frac#1#2{{#1 \over #2}}
%\headline={\rm  \ifodd\pageno  \RightHead  \else  \LeftHead  \fi}
%\def\RightHead{\centerline{
%Title
%}}
%\def\LeftHead{ \centerline{Doron Zeilberger}}
%end macros
\centerline
{\bf 
Estimating Many Constants With a Coin
}
\smallskip
\centerline
{\it Shalosh B. EKHAD and Doron ZEILBERGER}

\smallskip

\bigskip
{\bf Abstract:} We extend, and fully implement, in Maple, Jim Propp's charming way of estimating Pi via tossing a fair coin,
and compute many other constants, even going beyond the far more general set-up of Bruss and Paindaveine.

\bigskip

{\bf Maple packages}: This article is accompanied by two Maple packages {\tt Propp.txt}, and {\tt MultiProo.txt'}, available from

{\tt https://sites.math.rutgers.edu/\~{}zeilberg/tokhniot/Propp.txt } \quad ,

{\tt https://sites.math.rutgers.edu/\~{}zeilberg/tokhniot/MultiPropp.txt } \quad .

The front of this article

{\tt https://sites.math.rutgers.edu/\~{}zeilberg/mamarim/mamarimhtml/propp.html } \quad,

contains numerous output files, some of which will be referred to later.

{\bf Jim Propp's Way of Estimating Pi by Tossing a (fair) coin}

In  a delightful little note [P], Jim Propp described how to estimate Pi by tossing a coin. At each round, keep tossing a fair coin until the number of Heads exceeds the number of Tails,
then record the ratio:

{\it number of Heads divided by the number of Tosses} \quad .

So with probability $\frac{1}{2}$ you exit after one toss, and record $\frac{1}{1}=1$; with probability
$\frac{1}{8}$ you exit after $3$ tosses, recording $\frac{2}{3}$, etc. You do it for many rounds, say $10000$, and take the average of these recorded fractions.
What you get is a (not very good) {\it estimate} of $\pi/4$.

Jim must not have been aware of another nice paper [BP], that handles the much more general case of a {\it loaded} coin, $Pr(H)=p$, ($p \geq \frac{1}{2}$) and tossing the coin until the number of Heads exceeds the number of Tails
by $d+1$, the original case being $d=0$ and $p=\frac{1}{2}$. They derived, using human ingenuity, explicit expressions for the expected value of the ratio of the number of Heads to the total number of tosses. They also
had integral expressions for the variance.

The Maple package {\tt Propp.txt} goes even further, and shows how to efficiently compute averages of many other quantities at the termination, not just the above ratio, and in particular we will
present an alternative approach to [BP] using {\it algorithmic proof theory} [PWZ] and the {\it Zeilberger algorithm} [Z].

\vfill\eject

{\bf Jim Propp's original sum}

If the first time the number of Heads exceeded the number of Tails by $1$ was after tossing $k+1$ Heads and $k$ Tails, then the probability of that happening is
$C_k \cdot (\frac{1}{2})^k (\frac{1}{2})^{k+1} \, = \, C_k  (\frac{1}{2})^{2k+1}$,
where $C_k$ is the number of ways of tossing a coin $2k$ times, with $k$ Heads and $k$ Tails,
in such a way that the number of Heads never exceeds the number of Tails during the process. This number is famously the {\it Catalan number},
$$
C_k:=\frac{(2k)!}{k!(k+1)!} \quad .
$$
Since we are recoding the ratio of the number of Heads divided by the total number of tosses, we record $\frac{k}{2k+1}$. The expected value is
$$
\sum_{k=0}^{\infty} \, C_k (\frac{1}{2})^{2k+1} \cdot \frac{k}{2k+1} \quad, 
$$
and both Maple and Jim Propp agree that this equals $\frac{\pi}{4}$. Indeed, in Maple:

{\tt sum(1/(2*k+1)*(1/2)**(2*k+1)*(2*k)!/k!**2,k=0..infinity);}

outputs, after one nano-second, $\frac{\pi}{4}$.

{\bf Generalizing}

But what about a {\bf loaded coin}, with $Pr(Head)=p$ and $Pr(Tails)=1-p$?
And what about instead of recording $\#H/(\#H+ \#T)$ at the end of each run,
we record instead (for a {\bf fixed} $r$), 
$$
\frac{(\#H)_r}
{(\#H+ \#T)_r} \quad,
$$
where $(a)_r$ is the {\bf rising factorial}
$$
(a)_r \, := \, a(a+1) \cdots (a+r-1) \quad.
$$

Then the expected value of our `experiment', let's call it $A_r$ is
$$
A_{r} = \sum_{k=0}^{\infty} \frac{(2k)!^2 (k+r)! }{k!^2 (k+1)! (2k+r)!} \cdot p^{k+1}\, (1-p)^{k} \quad.
$$

Suppose that you want to compile a table of $A_r$ for $r$ from $1$ to $10000$. It would be very tedious to compute this infinite sum every time from scratch. In fact even for the
simple case $p=\frac{1}{2}$ and $r=2$, Maple can't evaluate it exactly, but returns it as a hypergeometric sum,
rather than the exact value $\frac{33 \pi}{64} \,- \, 1 $ \quad . We can compile such a table very efficiently, thanks to a third-order recurrence
satisfied by $A_r$, that is gotten immediately via the {\bf Zeilberger algorithm} [Z][PWZ].

\vfill\eject

The output file

{\tt https://sites.math.rutgers.edu/\~{}zeilberg/tokhniot/oPropp1.txt } \quad ,

that was obtained from the input file

{\tt https://sites.math.rutgers.edu/\~{}zeilberg/tokhniot/inPropp1.txt},

states this recurrence and lists the expected values for $r=1$ to $r=100$, for $p=\frac{3}{4}$, and compares them
to simulation estimates, running the `experiment' $10000$ times and averaging. It is fun to see the close agreement. For example with $r=74$, the
exact value is $0.75381265\dots$, while the simulation estimate is $0.7540493317\dots$. Not bad!

{\bf The Bruss-Paindaveiene generalization}

Let's go back to the case $r=1$, i.e. keeping track of $\#H/(\#H+ \#T)$, but {\bf now} you stop only when the number of Heads exceeds the number of Tails by $d+1$.
Now the expected value, let's call it $B_d$:
$$
B_d = \sum_{k=0}^{\infty} \, C(k,d)\, (1-p)^k \, p^{k+d+1} \frac{k+1}{2\,k+d+1} \quad,
$$
where $C(k,d)$ is the {\bf generalized Catalan number} (aka {\it ballot number})
$$
C(k,d)= \frac{(d+1) (2k+d)!}{k!(k+d+1)!} \quad.
$$

In a delightful paper [BP], F. Thomas Bruss and Dvy Paindaveine use sophisticated probability arguments to give an integral formula for $B_d$, and then express it
as an {\it explicit} but complicated infinite sum, that makes compiling a table for $B_d$, for say, $0 \leq d \leq 1000$, very time-consuming. Once again,
the Zeilberger algorithm comes to the rescue and finds a certain fourth-order linear recurrence.

This recurrene, along with the expected values for $0 \leq d \leq 70$, and  comparisons  to simulation outcomes can be found here:

{\tt https://sites.math.rutgers.edu/\~{}zeilberg/tokhniot/oPropp2.txt } \quad ,

{\bf Rolling a Die}

The Maple package

{\tt https://sites.math.rutgers.edu/\~{}zeilberg/tokhniot/MultiPropp.txt}

finds expected values and does simulate for the multi-dimensional analogs, but
is still very preliminary. You are welcome to extend it.

\vfill\eject

{\bf References}

[BP] F. Thomas Bruss and Dvy Paindaveine, {\it Win Rates at First-Passage Times for Biased Simple Random Walks}, 	arXiv:2512.21254 [math.PR]. \hfill\break
{\tt https://arxiv.org/abs/2512.21254}

[PWZ] Marko Petkovsek, Herbert S. Wilf, and Doron Zeilberger, {\it ``A=B''}, A.K. Peters.\hfill\break
{\tt https://sites.math.rutgers.edu/\~{}zeilberg/AeqB.pdf}

[P] Jim Propp, {\it Estimating $\pi$ with a Coin}, arXiv:2602.14487 [math.PR] . \hfill\break
{\tt https://arxiv.org/abs/2602.14487} \quad .

[Z] Doron Zeilberger, {\it The Method of Creative Telescoping}, J. Symbolic Computation {\bf 11} (1991) , 195-204. \hfill\break
{\tt https://sites.math.rutgers.edu/\~{}zeilberg/mamarimY/ct1991.pdf}

\bigskip
\hrule
\bigskip
Shalosh B. Ekhad and Doron Zeilberger, Department of Mathematics, Rutgers University (New Brunswick), Hill Center-Busch Campus, 110 Frelinghuysen
Rd., Piscataway, NJ 08854-8019, USA. {\tt ShaloshBEkhad@gmail.com, DoronZeil@gmail.com} \hfill\break
{\bf Aug. 20, 2026} 

\end